\documentclass[12pt,reqno]{article}

\usepackage[usenames]{color}
\usepackage{amssymb}
\usepackage{graphicx}
\usepackage{amscd}

\usepackage{amsthm}
\newtheorem{theorem}{Theorem}

\newtheorem{conjecture}[theorem]{Conjecture}

\theoremstyle{definition}

\newtheorem{example}[theorem]{Example}

\usepackage[colorlinks=true,
linkcolor=webgreen, filecolor=webbrown,
citecolor=webgreen]{hyperref}

\definecolor{webgreen}{rgb}{0,.5,0}
\definecolor{webbrown}{rgb}{.6,0,0}

\usepackage{color}
\usepackage{float}
\usepackage{graphics,amsmath,amssymb}
\usepackage{amsfonts}
\usepackage{latexsym}
\usepackage{epsf}

\newcommand{\seqnum}[1]{\href{http://oeis.org/#1}{\underline{#1}}}

\begin{document}

\begin{center}
\vskip 1cm{\LARGE\bf Integer sequences from elliptic curves} \vskip 1cm \large
Paul Barry\\
School of Science\\
South East Technological University\\
Ireland\\
\href{mailto:pbarry@wit.ie}{\tt pbarry@wit.ie}
\end{center}
\vskip .2 in

\begin{abstract} We indicate that given an integer coordinate point on an elliptic curve $$y^2+ a xy +by=x^3+cx^2+dx+e$$ we can identify an integer sequence whose Hankel transform is a Somos $4$ sequence, and whose Hankel determinants can be used to determine the coordinates of the multiples of this point. In reverse, given the coordinates of the multiples of an integer point on such an elliptic curve, we conjecture the form of a continued fraction generating function that expands to give a sequence with the above properties. \end{abstract}

\section{Introduction}
The interplay between elliptic curves, division polynomials, elliptic divisibility sequences, and Somos $4$ sequences \cite{GenCat, Chang, Hone, Poorten, Propp, Shipsey, Somos, Swart, Wang, Xin, Yura} has attracted much attention since the first paper by Ward \cite{Ward}. The work of Hone \cite{Hone} in particular has provided a uniform approach and solution, using the $\sigma$ function approach. The intention of this note is to show the fruitful interplay between integer points on an elliptic curve and integer sequences. These integer sequences can be defined starting from the equation of the curve; and the Hankel transforms \cite{Kratt1, Kratt2, Layman, Xin} of these sequences then provide Somos $4$ sequences. We remark that a classical link between Hankel determinants and elliptic curves (and hence Somos $4$ sequences following the results of Hone) is attributable to Kiepert, Frobenius and Stickelbeger \cite{Kiepert, Koecher, Onishi}. This says that
$$\frac{\sigma(nu)}{\sigma(u)^{n^2}}=\frac{1}{(-1)^{n-1}(1!2! \cdots (n-1)!)^2} \left|\begin{array}{rrrr}
\wp'(u) & \wp''(u) & \cdots & \wp^{(n-1)}(u)  \\
\wp''(u) & \wp'''(u) & \cdots & \wp^{(n)}(u)  \\
\vdots & \vdots & \vdots & \vdots  \\
\wp^{(n-1)}(u) & \wp^{(n)}(u) & \cdots & \wp^{(2n-3)}(u)  \\
\end{array}\right|.$$

\section{Example}
In order to motivate the form of the conjectures that are the main content of this note, we look at the case of the elliptic curve defined by
$$E:y^2+2xy+5y=x^3+4x^2+9x.$$
We solve this quadratic equation for $y$, to obtain
$$y = -\frac{\sqrt{4x^3+20x^2+56x+25}+2x+5}{2}, $$ which expands to give a sequence that begins
$$-5, -\frac{19}{5}, \frac{71}{125}, -\frac{2613}{3125},\ldots.$$
We have chosen this branch because the third term is positive. In order to work with integer values, we re-scale as follows.
$$\tilde{y} = -\frac{\sqrt{4(5^2x)^3+20(5^2x)^2+56(5^2x)+25}+2(5^2x)+5}{2(5)}, $$ to obtain a sequence that begins
$$-1,-19,71,-2613,78205,\ldots.$$
The sequence we seek will be tied to all elliptic curves with the same discriminant as the above curve; this means in practical terms that we start our sequence from the third term. Thus we consider the generating function
$$\frac{-\frac{\sqrt{4(5x)^3+20(5x)^2+56(5x)+25}+2(5x)+5}{10}-(-1-19x)}{x^2}=\frac{1+28x-\sqrt{1+56x+500x^2+2500x^3}}{2x^2}.$$ We next form the generating function
$$\frac{1}{1-x-x^2\left(\frac{1+28x-\sqrt{1+56x+500x^2+2500x^3}}{2x^2}\right)},$$ or
$$\frac{1}{1-30x+\sqrt{1+56x+500x^2+2500x^3}}.$$
We now \emph{revert} this generating function to obtain
$$u(x)=\frac{1+30x+100x^2-\sqrt{1+60x+1100x^2+3500x^3-62500x^4}}{1250x^3}.$$
We finally form the generating function
$$\frac{1}{1-x-x^2 u(x)}=\frac{1250x}{\sqrt{1+60x+1100x^2+3500x^3-62500x^4}-1350x^2+1220x-1}.$$
This is the generating function of the integer sequence that we seek. In this case, it expands to give a sequence $a_n$ that begins $$1,1,2,2,-67,2688, -73696, 1856194, \ldots.$$ The Hankel transform of this sequence is the sequence $h_n=|a_{i+j}|_{0 \le i,j \le n}$ of Hankel determinants. We are also interested in the modified Hankel determinants $h_n^*=|\tilde{a}_{i,j}|_{0 \le i, j \le n}$ where
\begin{align*}\tilde{a}_{i,j}=\begin{cases} a_{i+j} &\text{if $i<n$}, \\
a_{i+j+1} &\text{otherwise}. \end{cases}\end{align*}
The re-scaled Hankel transform $\tilde{h}_n=\frac{h_n}{5^{n^2-2n}}$ begins
$$1,5,-71,-13065, -1275214, 2876558965,\ldots.$$
This sequence satisfies the identity
$$\tilde{h}_n = \frac{ 25 \tilde{h}_{n-1} \tilde{h}_{n-3} + 71 \tilde{h}_{n-2}^2}{\tilde{h}_{n-4}}.$$
This means that it is a $(25, 71)$ Somos-$4$ sequence. (Note that the third term in the expansion of $y$, namely $\frac{71}{5^3}$, already contains these parameters).

Now the \emph{division polynomial sequence} of the elliptic curve $E$ begins
$$0,1,-5,-71,13065, -1275214, -2876558965,\ldots.$$
Thus, up to sign, we have arrived at the elliptic divisibility of the curve using generating functions and Hankel transforms, starting with the solution of the quadratic equation in $y$. Significantly, we have found an integer sequence. Indeed, more can be found at this stage. The point $(0,0)$ lies on the above elliptic curve $E$, and the above sequence $0,1,-5,-7,\ldots$ is attached to this point. The point $(0,-5)$ is the other integer point on this curve, and the sequence $0,1,5,-71,-1305,\ldots$ corresponds to this. In both cases, the sequence is dependent on the coordinates of the integer multiples of the original point. We now find the following. The $x$ coordinates of the multiples of $(0,0)$ on $E$ are given by
\begin{align*}x_n=\begin{cases} 0, & \text{if $n=0$},\\
-\frac{1}{5^2} \frac{h_{n-1} h{_n+1}}{h_n^2},&\text{otherwise}.\end{cases}\end{align*}
Similarly, the $y$ coordinates of the integer multiples of $(0,0)$  on $E$ are given by
\begin{align*}y_n=\begin{cases} 0, &\text{if $n=0$},\\
-\frac{1}{5^3} \frac{h_{n-1} h{_n+1}}{h_n^2}\left(\frac{h_{n+1}^*}{h_{n+1}}-\frac{h_n^*}{h_n}+9+1\right), &\text{otherwise}.\end{cases}\end{align*} Note that the numbers $5$ and $9$ here are the parameters appearing in the defining equation
$$E: y^2+2xy+5y=x^3+4x^2+9x.$$

The generating function of the sequence $a_n$ can be expressed as a Jacobi continued fraction \cite{CFT, Wall}
$$\cfrac{1}{1-\alpha_0 x-\cfrac{\beta_1 x^2}{1-\alpha_1 x - \cfrac{\beta_2 x^2}{1- \alpha_2 x-\cdots}}},$$
where
$$\alpha_n=\frac{h_n^*}{h_n}-\frac{h_{n-1}^*}{h_{n-1}}+0^n,$$
and
\begin{align*}\beta_n =\begin{cases} 1, &\text{if $n=0$},\\
 \frac{h_{n-1}h_{n+1}}{h_n^2}, &\text{otherwise}.\end{cases}\end{align*}

In reverse, if we start with the coefficients $(x_n,y_n)=n(0,0)$ on $E$ then we can retrieve the sequence $a_n$ through its generating function in continued fraction form as follows. For this, we set
\begin{align*}\alpha_n = \begin{cases} 1, &\text{if $n=0$},\\
 -1,  &\text{if $n=1$},\\
 \frac{5 y_{n-1}}{x_{n-1}}-(9+1), &\text{otherwise},\end{cases}\end{align*}
 and
\begin{align*}\beta_{n+1} = \begin{cases} 1, &\text{if $n=0$},\\
 -5^2 x_n, &\text{otherwise}.\end{cases}\end{align*}

\begin{example} The elliptic curve $E:y^2+2xy+5y=x^3+4x^2+9x$ has discriminant equal to $-38091$. The elliptic curve $\tilde{E}:y^2+5y=x^3+5x^2+14x$ also has discriminant $-38091$. Solving for $y$ for $\tilde{E}$, we obtain
$$y=-\frac{5+\sqrt{4x^3+20x^2+56x+25}}{2},$$ which expands to give the sequence that begins
$$-5, -\frac{14}{5}, \frac{71}{125}, -\frac{2613}{3125},\ldots.$$ We see that after the second term, the sequences for $E$ and $\tilde{E}$ coincide. The curve with Cremona label $38091$a$1$ \cite{LMFDB} is the curve
$$y^2+y=x^3-x^2+6x-10$$ with integral points $(2,2)$ and $(2,-3)$. For the point  $(2,2)$, we obtain the division polynomial sequence $0,1,-5,-71,13065,\ldots$.
\end{example}

\section{A general case}
We now consider the case of the elliptic curve
$$E=E(a,b,c,d): y^2+axy+by=x^3+cx^2+dx.$$ The point $(0,0)$ is clearly on this curve, and we shall use this as the reference point. The discriminant of this curve,
\begin{scriptsize}
$$dba^5 + (-b^2c + d^2)a^4 + (8dbc + b^3)a^3 + (-8b^2c^2 + 8d^2c - 30db^2)a^2 + (16dbc^2 + 36b^3c - 96d^2b)a + (-16b^2c^3 + 16d^2c^2 + 72d*b^2c + (-27b^4 - 64d^3))$$ \end{scriptsize} is assumed to be nonzero. Following the steps outlined above (where we use $b$ as the scaling factor in place of $5$ above), we arrive at the generating function
$$g(x)=\frac{2b^4x}{\sqrt{Ax^4+Bx^3+Cx^2+Dx+1}+Fx^2+Gx-1},$$ where
\begin{align*}
A&=a^2b^2(d + 1)^2 - 2ab(2b^4 + b^2c(d + 1) - (d + 1)^3) + b^4(c^2 - 4(2d + 1)) - 2b^2c(d + 1)^2 + (d + 1)^4,\\
B&=2(a^2b^2(d + 1) + ab(3(d + 1)^2 - b^2c) - 2(b^4 + b^2c(d + 1) - (d + 1)^3)),\\
C&=a^2b^2 + 6ab(d + 1) - 2(b^2c - 3(d + 1)^2),\\
D&=2(ab + 2(d + 1)),\\
F&=- (ab(d + 1) + 2b^4 - b^2c + (d + 1)^2),\\
G&=2(b^4 - d - 1) - ab.\end{align*}
This expands to give a sequence $a_n$ which begins
$$1, 1, 2, 2, - abd + b^2c - d^2 + 4, a^2b^2d + ab(d(3d + 1) - b^2c) + b^4 - b^2c(2·d + 1) + 2d^3 + d^2 + 4,\ldots.$$
This then has a Hankel transform $h_n$ such that $\frac{h_n}{b^{n^2-2n}}$ begins
$$1, b, - abd + b^2c - d^2, - b(a^2b^2d + ab(3d^2 - b^2c) + b^4 - 2·b^2cd + 2d^3), \ldots.$$
We then have the following  conjecture.
\begin{conjecture} The sequence $\frac{h_n}{b^{n^2-2n}}$ is a $(b^2, abd - b^2c + d^2)$ Somos-$4$ sequence which coincides with the division polynomial sequence of $E$.
\end{conjecture}
We also have the following conjecture.
\begin{conjecture}
The $x$ coordinates of the multiples of $(0,0)$ on $E$ are given by
\begin{align*}x_n=\begin{cases} 0, & \text{if $n=0$},\\
-\frac{1}{b^2} \frac{h_{n-1} h{_n+1}}{h_n^2}, & \text{otherwise}.\end{cases}\end{align*}
Similarly, the $y$ coordinates of the integer multiples of $(0,0)$  on $E$ are given by
\begin{align*}y_n=\begin{cases} 0, &\text{if $n=0$},\\
-\frac{1}{b^3} \frac{h_{n-1} h{_n+1}}{h_n^2}\left(\frac{h_{n+1}^*}{h_{n+1}}-\frac{h_n^*}{h_n}+d+1\right), &\text{otherwise}.\end{cases}\end{align*}
\end{conjecture}
In the reverse direction, we assume that $(x_n, y_n)$ are the coordinates of $n(0,0)$ on $E=E(a,b,c,d)$. Then  we have the following conjecture.
\begin{conjecture} The generating sequence $g(x)$ can be expressed as the Jacobi continued fraction
$$g(x)=\cfrac{1}{1-\alpha_0 x - \cfrac{\beta_1 x^2}{1-\alpha_1 x - \cfrac{\beta_2 x^2}{1-\alpha_2 x- \cdots}}},$$ where
\begin{align*}\alpha_n = \begin{cases} 1, &\text{if $n=0$},\\
 -1, &\text{if $n=1$},\\
 \frac{b y_{n-1}}{x_{n-1}}-(d+1), &\text{otherwise},\end{cases}\end{align*}
 and
\begin{align*}\beta_{n+1} = \begin{cases} 1, &\text{if $n=0$},\\
 -b^2 x_n, &\text{otherwise}.\end{cases}\end{align*}
 \end{conjecture}

\section{Singular cubic curves}
For completeness we apply our formalism to two singular (non-elliptic) cubic curves.
\begin{example} We consider the singular cubic curve
$$y^2+xy+y=x^3-2x^2$$ which has a node at the point $(1,-1)$. Solving for $y$, we obtain
$$y=\frac{(x-1)\sqrt{1-4x}-x-1}{2}$$ which expands to give the sequence
$$-1, -1, 2, -3, 7, -19, 56, -174, 561, \ldots.$$
Note that the sequence $C_n+C_{n+1}$ begins $2,3,7,19,\ldots$ where $C_n=\frac{1}{n+1}\binom{2n}{n}$ is the $n$-th Catalan number \seqnum{A000108}.
We now form
$$\frac{\frac{(x-1)\sqrt{1-4x}-x-1}{2}+1+x}{x^2}=\frac{1+x-(1-x)\sqrt{1+4x}}{2x^2}.$$
Forming the generating function
$$\frac{1}{1-x-x^2\left(\frac{1+x-(1-x)\sqrt{1+4x}}{2x^2}\right)}=\frac{2}{1-3x+(1-x)\sqrt{1+4x}}$$ and
reverting, we obtain the generating function
$$\tilde{g(x)}=\frac{1+3x+4x^2-(1+x)\sqrt{1+4x+8x^2}}{2x^3}.$$
We finally form the generating function $\frac{1}{1-x-x^2 \tilde{g}(x)}$ to obtain
$$g(x)=\frac{2x}{(1+x)\sqrt{1+4x+8x^2}-6x^2-x-1}.$$
This expands to give a sequence that begins
$$1, 1, 2, 2, 2, 9, -7, 25, 19, -125, 474,\ldots.$$
The Hankel transform of this sequence is then given by the signed Fibonacci numbers $(-1)^{\binom{n}{2}} F_{n+1}$ (\seqnum{A000045}):

$$1, 1, -2, -3, 5, 8,-13,-21,34,\ldots.$$
This is a $(1,2)$ Somos-$4$ sequence.
\end{example}
\begin{example} We consider the singular cubic curve
$$y^2+2y=x^3-x^2-x,$$ which again has a node at $(1,-1)$. We find that
$$g(x)=\frac{32x}{(1-2x)\sqrt{1+4x+20x^2}-36x^2+32x-1}.$$ This expands to give a sequence that begins
$$1, 1, 2, 2, -1, 15, 8, -152, 493, 541, -8898,\ldots.$$
With $h_n$ signifying the Hankel transform of this sequence, we find that
$\frac{h_n}{2^{n^2-2n}}$ is the sequence that begins
$$1, 2, -5, -12, 29, 70,\ldots.$$ This is the signed Pell sequence $(-1)^{\binom{n}{2}}P_{n+1}$ (\seqnum{A000129}).
This is a $(4, 5)$ Somos-$4$ sequence.
\end{example}
\begin{example} The cubic curve $y^2=x^3+x^2$ has a node at $(0,0)$. The methods used so far in this note do not work in this instance.
\end{example}

\section{Riordan arrays and closed expressions}
Using the theory of Riordan arrays, we can use the so-called ``fundamental theorem of Riordan arrays'' \cite{book, SGWW} to express the generating function $g(x)$ in a form that allows us to find a closed form expression for the elements $a_n$ of the expansion of $g(x)$. For this, we let
\begin{align*}
\alpha&=ab-2(b^4-d-1)\\
\beta&=ab(d+1)+2b^4-b^2c+(d+1)^2\\
\gamma&=ab(d+2)+b^4-b^2c+d^2+4d+2\\
\delta&=abd+2b^4-b^2c+d^2-2\\
\epsilon&=ab-b^4+2d+1.\end{align*}
Then we can express $g(x)$ as
$$\frac{b^4x}{1+\alpha x + \beta x^2}c\left(\frac{b^4x(\gamma x^3-\delta x^2- \epsilon x-1)}{(1+\alpha x + \beta x^2)^2}\right)+\frac{1+\alpha x + \beta x^2}{1+\epsilon x + \delta x^2- \gamma x^3},$$
where $c(x)=\frac{1-\sqrt{1-4x}}{2x}$ is the generating function of the Catalan numbers. We let
$$S(r,\delta, \epsilon, \gamma)=\sum_{i=0}^r \sum_{j=0}^{r-i}\binom{i}{j}\binom{j}{r-i-j}\epsilon^{i-j}(-\gamma)^{r-i-j} \delta^{2j+i-r},$$
which is the expansion of $\frac{1}{1+\epsilon x + \delta x^2- \gamma x^3}$.
Then we have
$$a_n=b^4 \sum_{k=0}^{n-1}b^{4k}(-1)^k \sum_{j=0}^k \binom{k}{j}\sum_{l=0}^j\binom{j}{l} \epsilon^{j-l}\sum_{r=0}^l \binom{l}{r}(-\gamma)^r \delta^{l-r}\sum_{i=0}^{n-k-j-r-l}\binom{2k+i}{i}$$
$$\quad\quad\quad \cdot \binom{i}{n-k-j-l-r-i-1}(-1)^i \beta^{n-k-j-l-r-i-1}\alpha^{2i+r+l+j+k-n+1}$$
$$\quad + S(n,\delta,\epsilon,\gamma)+\alpha S(n-1,\delta,\epsilon,\gamma)+\beta S(n-2,\delta,\epsilon,\gamma).$$

\section{A note on the integer sequences}
By construction, all sequences with generating function $g(x;a,b,c,d)$ will start $1,1,2,2,\ldots$. This does not mean that only such sequences will have the required Hankel transform, since many sequences may have the same Hankel transform. In particular, given a sequence $a_n$ with generating function $g(x)$, then the sequence with generating function $\frac{1}{1-rx} g\left(\frac{x}{1-rx}\right)$ (binomial transform), and the sequence with generating function $\frac{g(x)}{1-r x g(x)}$ (INVERT transform), or any combination of these transforms, will have the same Hankel transform. Thus many sequences appear in the literature (and most notably in the On-Line Encyclopedia of Integer Sequences \cite{SL1, SL2}) with a different form.
\begin{example} We consider sequence \seqnum{A178072} which begins
$$1, 0, -1, -1, -1, -1, 1, 8, 23, 45, 55, -14,\ldots,$$ and which has its generating function given by
$$g_0(x)=\frac{2}{1+2x+x^2+\sqrt{1-4x+6x^2+x^4}}.$$
The Hankel transform of this sequence begins
$$1, -1, 1, 2, -1, -3, -5, 7, -4, -23, 29, 59,\ldots,$$ which is to be compared with \seqnum{A006769}, the elliptic divisibility sequence for the elliptic curve $y^2+y=x^3-x$. This elliptic divisibility sequence begins
$$0, 1, 1, -1, 1, 2, -1, -3, -5, 7, -4, -23, 29, 59, 129, -314,\ldots.$$ The theory of Hankel transforms now tell us that the sequence with generating function
$$g_1(x)=\frac{1}{1-x-x^2 g_0(x)}=\frac{4-5x-3x^3-x \sqrt{1-4x+6x^2+x^4}}{2(2-5x+2x^2+2x^4)}$$ will have a Hankel transform that begins
$$1, 1, -1, 1, 2, -1, -3, -5, 7, -4, -23, 29, 59, 129, -314,\ldots.$$
This sequence begins
$$1, 1, 2, 3, 4, 5, 5, 3, -1, -3, 12, 79, 253, 565, 858, \ldots.$$
Our theory tells us that the sequence with generating function
$$g(x)=\frac{1-2x+2x^2+\sqrt{1-4x^3+4x^4}}{2(1-x)^2}$$ will have an equal Hankel transform. This sequence begins $$1, 1, 2, 2, 3, 4, 4, 6, 7, 6, 11,\ldots.$$
In order to compare these two generating functions ($g(x)$ and $g_1(x)$) we note that if a generating function has a Jacobi continued fraction expression
$$\cfrac{1}{1-\alpha_0 x-\cfrac{\beta_1 x^2}{1-\alpha_1 x - \cfrac{\beta_2 x^2}{1- \alpha_2 x-\cdots}}},$$ then
the corresponding Hankel transform depends only on the coefficients $\beta_n$, and thus by varying the coefficients $\alpha_n$, it is seen that many sequences can have the same Hankel transform. In this case, if we let $\alpha_n$ be associated to $g(x)$, and $\alpha_n^{(1)}$ be associated to $g_1(x)$, then we have
that $\alpha_n + \alpha_n^{(1)}$ is the sequence that begins
$$2, -1, 1, 1, 1, 1, 1, 1, 1, 1,\ldots.$$ This shows the relationship between the $\alpha$-sequences. The $\beta$-sequences are of course identical.
\end{example}
\begin{example} We consider the elliptic curve $$E_1: y^2+xy=x^3-2x+1$$ which has an integer point at $(1,-1)$. Translating this to $(0,0)$, we obtain the curve
$$E: y^2+xy-y=x^3+3x^2+2x.$$
For this curve, we find that
$$g(x)=\frac{1+3x+5x^2-\sqrt{1+10x+31x^2+26x^3-7x^4}}{2(1+3x-x^2-8x^3)}.$$
This expands to give the sequence $a_n$ that begins
$$1, 1, 2, 2, 5, 1, 24, -53, 278, -1048, 4442,\ldots.$$
The Hankel transform of this sequence is the $(1,-1)$ Somos-$4$ sequence (see \seqnum{A178079}) that begins
$$1, 1, 1, 2, 1, -3,-7,-8,-25,\ldots.$$
The sequence $a_n$ with generating function $g(x)$ is closely related to the known sequence \seqnum{A178078}, which has generating function
$$g_1(x)=\frac{1+3x-x^2-\sqrt{1-6x+7x^2+2x^3+x^4}}{2x(3-2x^2)}.$$ The Hankel transform of this sequence is given by
$$1, 1, 2, 1, -3,-7,-8,-25,\ldots.$$
The same will be true for the generating function $\tilde{g}(x)$ where
$$g(x)=\frac{1}{1-x-x^2 \tilde{g}(x)}.$$ Here, we have that 
$$\tilde{g}(x)=\frac{1+5x+3x^2-\sqrt{1+10x+31x^2+26x^3-7x^4}}{2x^3}.$$
We have the following relationship: $g_1(x)$ is the INVERT$(-3)$ transform of the $4^{\text{th}}$ binomial transform of $\tilde{g}(x)$.

$$g_1(x)=\frac{\frac{1}{1-4x} \tilde{g}\left(\frac{x}{1-4x}\right)}{1+3x\frac{1}{1-4x} \tilde{g}\left(\frac{x}{1-4x}\right)}.$$

For completeness, we point out a technical issue that arises in the process, starting with the equation $y^2+xy-y=x^3+3x^2+2x$, that leads to the generating function $g(x)$ in this example. Solving for $y$, and dropping the first two terms of the expansion, we obtain the term
$$\frac{\frac{1-x+\sqrt{1+6x+13x^2+4x^3}}{2}-1-x}{x^2}=\frac{\sqrt{1+6x+13x^2+4x^3}-3x-1}{2x^2}.$$
We next calculate
$$\frac{1}{1-x+x^2\left(\frac{\sqrt{1+6x+13x^2+4x^3}-3x-1}{2x^2}\right)}=\frac{2}{1-5x+\sqrt{1+6x+13x^2+4x^3}}.$$ The choice of the $'+'$ sign has given us the desired form of expression at this stage.
\end{example}

\section{Comment on reversion}
If a power series $f(x)=a_0+a_1x+a_2x^2+\cdots$ is such that $a_0=0, a_1 \ne 0$, then we can find its reversion or its compositional inverse $\bar{f}(x)=v(x)$, where $v$ is the solution of the equation $f(v)=x$ such that $v(0)=0$. For a generating function $g(x)=a_0+a_1x+ a_2 x^2+\cdots$ where $a_0 \ne 0$, we define its reversion to be $\frac{1}{x} \overline{xg}$. Using Lagrange inversion, the coefficients of such a reversion are given by
$$\frac{1}{n+1} [x^n] \frac{1}{g(x)^{n+1}},$$ that is, the reversion of $g(x)$ is given by
$\sum_{n=0}^{\infty}\frac{1}{n+1} [t^n] \frac{1}{g(t)^{n+1}} x^n$.

\section{A fuller picture}
So far, we have concentrated on the full process that leads to the desired generating function $g(x)$. It is nevertheless of interest to examine an intermediate result, and to place this in a broader context. For this, we take the special case of the elliptic curve $y^2+y=x^3-x$ (LMFDB label $37.a1$). Solving for $y$, we find
$$y=-\frac{1+\sqrt{1-4x+4x^2}}{2},$$ which expands to give
$$-1, 1, 1, 1, 3, 8, 23, 68, 207,\ldots.$$ We thus look at the generating function
$$\frac{-\frac{1+\sqrt{1-4x+4x^2}}{2}+1-x}{x^2}=\frac{1-2x-\sqrt{1-4x+4x^3}}{2x^2},$$ which expands to \seqnum{A056010} (and \seqnum{A025262}$(n+1)$),
$$1, 1, 3, 8, 23, 68, 207, 644, 2040, 6558, 21343, 70186, 232864,\ldots.$$ The Hankel transform of this sequence begins
$$1, 2, 3, 7, 23, 59, 314, 1529, 8209, 83313, 620297,\ldots.$$ This coincides with the Somos-$4$ sequence \seqnum{A006720}$(n+3)$.
The next step in the process is to form the generating function
$$\frac{1}{1-x-x^2\left(\frac{1-2x-\sqrt{1-4x+4x^3}}{2x^2}\right)}=\frac{2}{1+\sqrt{1-4x+4x^3}},$$ which expands to give the sequence \seqnum{A157003},
$$1, 1, 2, 4, 10, 27, 78, 234, 722, 2274, 7280,\ldots.$$ The Hankel transform of this sequence is the Somos-$4$ sequence \seqnum{A006720}$(n+2)$,
$$1, 1, 2, 3, 7, 23, 59, 314, 1529, 8209, 83313, 620297,\ldots.$$
We now revert the generating function $\frac{2}{1+\sqrt{1-4x+4x^3}}$ to obtain the generating function
$\frac{1-\sqrt{1-4x^3+4x^4}}{2x^3}$. This expands to give a sequence
$$1, -1, 0, 1, -2, 1, 2, -6, 6, 3, -20, 30, -6, -65, \ldots.$$ The Hankel transform in this case begins
$$1, -1, 1, 2, -1, -3, -5, 7, -4, -23, 29,\ldots,$$ or \seqnum{A006769}$(n+2)$, where \seqnum{A006769} is the elliptic divisibility sequence of the elliptic curve $ y^2 + y = x^3 - x$.

Finally, we form the generating function
$$g(x)=\frac{1}{1-x-x^2\left(\frac{1-\sqrt{1-4x^3+4x^4}}{2x^3}\right)}=\frac{1-2x+\sqrt{1-4x^3+4x^4}}{2(1-x)^2}.$$ This expands to give the sequence
$$1, 1, 2, 2, 3, 4, 4, 6, 7, 6, 11, 10, 6, 22, 8, 0, \ldots.$$ This sequence then has a Hankel transform
$$1, 1, -1, 1, 2, -1, -3, -5, 7, -4, -23, \ldots.$$ Apart from an initial term $0$, this is the elliptic divisibility sequence of $y^2+y=x^3-x$. The bisection of this sequence begins
$$1, -1, 2, -3, 7, -23, 59, -314, 1529, -8209, 83313,\ldots.$$ This is an alternating sign version of the Hankel transform of \seqnum{A157003} above.

\section{Conclusions} Given an integer point $(x,y)$ on an elliptic curve
$$E: y^2+ axy+by=x^3+cx^2+dx+e$$ where $a,b,c,d \in \mathbb{Z}$, we may translate the curve by $(x,y)$ to obtain a curve with equation
$$E': y^2+a' xy + b' y = x^3+c'x^2+d'x,$$  where the point of interest is now $(0,0)$. Using our formalism, we derive a generating function for a sequence whose scaled Hankel transform $\frac{h_n}{b^{2n-n}}$ coincides with the division polynomial of the curve. This sequence will be a Somos $4$ sequence. Using the Hankel parameters, we can also find the coordinates of the multiples of the point $(0,0)$ on $E'$. In the opposite direction, if we know the coordinates of the multiples of $(0,0)$ on such a curve, then we can use these to construct the Jacobi coefficients of a generating function that coincides with the generating function constructed from the equation of the curve. A critical role is played by the process of reversion of generating functions, though the reasons why this is so are still mysterious.

\bigskip
\hrule

\noindent 2010 {\it Mathematics Subject Classification}:
Primary 11B83; Secondary 14H52, 11B39,30B70
\noindent \emph{Keywords:} Elliptic curve, Somos sequence, division polynomial, elliptic divisibility sequence, Hankel determinant, integer sequence, Jacobi continued fraction, Riordan array.

\bigskip
\hrule
\bigskip
\noindent (Concerned with sequences
\seqnum{A000045},
\seqnum{A000108},
\seqnum{A000129},
\seqnum{A006720},
\seqnum{A006769},
\seqnum{A025262},
\seqnum{A056010},
\seqnum{A157003},
\seqnum{A178072},
\seqnum{A178078} and
\seqnum{A178079}).


\begin{thebibliography}{99}

\bibitem{book} P. Barry, \emph{Riordan Arrays: a Primer}, Logic Press, 2017.

\bibitem{GenCat} P. Barry, Generalized Catalan numbers, Hankel transforms and Somos-$4$ sequences, \emph{J. Integer Seq.} \textbf{13} (2010), \href{https://cs.uwaterloo.ca/journals/JIS/VOL13/Barry1/barry95r.html} {Article 10.7.2}.

\bibitem{CFT} P. Barry, Continued fractions and
transformations of integer sequences, \emph{J. Integer Seq.}, \textbf{12} (2009), \href{https://cs.uwaterloo.ca/journals/JIS/VOL12/Barry3/barry93.html}{Article 09.7.6}.

\bibitem{Chang} X.-K. Chang and X.-B. Hu, A conjecture based on Somos-$4$ sequence and its extension, \emph{Linear Algebra Appl.} \textbf{436} (2012), 4285--4295.

\bibitem{Hone} A. N. W. Hone, Elliptic curves and quadratic recurrence sequences, \emph{Bull. Lond. Math. Soc.} \textbf{37} (2005), 161--171.

\bibitem{Kiepert} L. Kiepert, Wirkliche Ausfurhrung der ganzzahligen Multiplikation der elliptischen Funktionen, \emph{J. reine angew. Math.} \text{76} (1873), 21--33.

\bibitem{Koecher} M. Koecher and A. Krieg, \emph{Elliptische Funktionen und Modulformen}, Springer.

\bibitem{Kratt1} C. Krattenthaler, Advanced determinant
    calculus, \emph{S\'eminaire Lotharingien Combin.} \textbf{42} (1999), Article B42q., available electronically at
    \texttt{https://arxiv.org/abs/math/0503507},
    2023.

\bibitem{Kratt2} C. Krattenthaler, Advanced determinant
    calculus: A complement, {\it Linear Algebra
    Appl.} \textbf{411} (2005), 68–-166.

\bibitem{Layman} J. W. Layman, The Hankel transform and some of its properties, \emph{J. Integer Seq.}, \textbf{4} (2001),
\href{https://www.cs.uwaterloo.ca/journals/JIS/VOL4/LAYMAN/hankel.html} {Article 01.1.5}.

\bibitem{LMFDB} The LMFDB Collaboration, The $L$-functions and modular forms database, \url{https://www.lmfdb.org}, 2023.


\bibitem{Onishi} Y. Onishi, Determinant expressions for hyperelliptic functions, \emph{Proc. Edinb. Math. Soc.} \textbf{48} (2005),  705--742

\bibitem{Poorten} A. van der Poorten and C. Swart, Recurrence relations for elliptic sequences: every Somos $4$ is a Somos $k$, \emph{Bull. Lond. Math. Soc.}, \textbf{38} (2006), 546--554.

\bibitem{Propp} J. Propp, The Somos sequence site, \url{http://jamespropp.org/somos.html}.

\bibitem{SGWW} L. W. Shapiro, S. Getu, W-J. Woan, and L.C. Woodson,
The Riordan group, \emph{Discr. Appl. Math.}, \textbf{34} (1991),
 229--239.

\bibitem{Shipsey} R. Shipsey, \emph{Elliptic Divisibility Sequences}, PhD Thesis, Goldsmiths, University of London, 2001.


\bibitem{SL1} N. J. A.~Sloane, \emph{The
On-Line Encyclopedia of Integer Sequences}. Published electronically
at \texttt{http://oeis.org}, 2023.

\bibitem{SL2} N. J. A.~Sloane, The On-Line Encyclopedia of Integer
Sequences, \emph{Notices Amer. Math. Soc.}, \textbf{50} (2003),  912--915.

\bibitem{Somos} M. Somos, \url{https://grail.cba.csuohio.edu/~somos/math.html}.


\bibitem{Swart} C. Swart, \emph{Elliptic Curves and Related Sequences}, PhD Thesis, Royal Holloway and Bedford New College, University of London, 2003

\bibitem{Wall} H. S. Wall, \emph{Analytic Theory of Continued Fractions}, AMS, 2000.

\bibitem{Xin} G. Xin, Proof of the Somos-$4$ Hankel determinants conjecture, \emph{Adv. in Appl. Math.} \textbf{42} (2009), 152--156.

\bibitem{Wang} Y. Wang and Z. Zhang, Proof of four $(\alpha, \beta)$ Somos $4$ Hankel determinants conjectures of Barry, available at \url{https://arxiv.org/abs/2305.05995}, 2023.

\bibitem{Yura} F. Yura, Hankel determinant solution for elliptic sequence, \emph{Linear Algebra Appl.}, \textbf{484} (2015), 27--45.

\bibitem{Ward} M. Ward, Memoir on elliptic divisibility sequences, \emph{Amer. J. Math.} \textbf{70} (1948), 31--74.


\end{thebibliography}
\end{document}